\documentclass{gtart}

\def\ifplaintex{\expandafter\ifx\csname documentclass\endcsname\relax}

\def\gtp{{\mathsurround=0pt\it $\cal G\mskip-2mu$eometry \&\ 
$\cal T\!\!$opology $\cal P\!$ublications}}  

\def\recd{{\small Received:\qua\receiveddate\ifx\reviseddate\relax
\else\qquad Revised:\qua\reviseddate\fi\par}} 

\def\lognumber#1{\def\thelognumber{#1}}
\def\volumenumber#1{\def\thevolumenumber{#1}}
\def\volumeyear#1{\def\thevolumeyear{#1}}
\def\papernumber#1{\def\thepapernumber{#1}}
\def\pagenumbers#1#2{\def\startpage{#1}\def\finishpage{#2}}
\def\published#1{\def\publishdate{#1}}

\def\received#1{\def\receiveddate{#1}}
\def\revised#1{\def\reviseddate{#1}}
\def\accepted#1{\def\accepteddate{#1}}

\def\asciiaddress#1{\def\theasciiaddress{#1}}

\long\def\asciiabstract#1{\long\def\theasciiabstract{#1}}

\let\\\par\let\thelognumber\relax\let\thevolumenumber\relax
\let\thepapernumber\relax\let\thevolumeyear\relax\let\startpage\relax
\let\finishpage\relax\let\publishdate\relax\let\receiveddate\relax
\let\reviseddate\relax\let\accepteddate\relax\let\theasciititle\relax
\let\theasciiauthors\relax\let\theasciiaddress\relax
\let\theasciiabstract\relax

\let\theasciiemail\relax

\ifplaintex
\font\logobig=cmssbx10 scaled 3836
\font\logomed=cmssbx10 scaled 2557
\else
\font\logobig=cmssbx10 scaled 4200
\font\logomed=cmssbx10 scaled 2800
\fi

\long\def\makeagttitle{   
\count0=\startpage
\agt\hfill      
\hbox to 45truept{\vbox to 0pt{\vglue -13truept{\logomed A\kern -.37em{\logobig 
T}\kern -.38em G}\vss}\hss}
\break
{\small Volume \thevolumenumber\ (\thevolumeyear)
\startpage--\finishpage\nl
Published: \publishdate}

\vglue .25truein

{\parskip=0pt\leftskip 0pt plus
1fil\def\\{\par\smallskip}{\Large\bf\thetitle}\par\medskip} \vglue
0.05truein

{\parskip=0pt\leftskip 0pt plus 1fil\def\\{\par}{\sc\theauthors}
\par\medskip}%
 
\vglue 0.03truein 

{\small\leftskip 25truept\rightskip 25truept{\bf Abstract}\stdspace\theabstract

{\bf AMS Classification}\stdspace\theprimaryclass
\ifx\thesecondaryclass\relax\else; \thesecondaryclass\fi\par
{\bf Keywords}\stdspace \thekeywords\par}\vglue 7truept

}   

\ifplaintex
\font\phead=cmsl9 scaled 950
\font\pnum=cmbx10 scaled 913
\font\pfoot=cmsl9 scaled 950
\headline{\vbox to 0pt{\vskip -4.5mm\line{\small\phead\ifnum
\count0=\startpage ISSN 1472-2739 (on-line) 1472-2747 (printed)
\hfill {\pnum\folio}\else\ifodd\count0\def\\{ }%
\ifx\theshorttitle\relax\thetitle\else\theshorttitle\fi\hfill{\pnum\folio}
\else\def\\{ and }{\pnum\folio}\hfill\ifx\theshortauthors\relax\theauthors
\else\theshortauthors\fi\fi\fi}\vss}}
\footline{\vbox to 0pt{\vglue 0mm\line{\small\pfoot\ifnum\count0=\startpage
\copyright\ \gtp\hfill\else
\agt, Volume \thevolumenumber\ (\thevolumeyear)\hfill\fi}\vss}}
\else
\font=cmsl9 scaled 1050
\font\lnum=cmbx10 
\font=cmsl9 scaled 1050
\makeatletter
\def\@oddhead{{\small\lhead\ifnum\count0=\startpage ISSN 1472-2739 
(on-line) 1472-2747 (printed)\hfill {\lnum\number\count0}\else\ifodd\count0
\def\\{ }\ifx\theshorttitle\relax \thetitle \else\theshorttitle\fi\hfill
{\lnum\number\count0}\else\def\\{ and }{\lnum\number\count0}
\hfill\ifx\theshortauthors\relax 
\theauthors\else\theshortauthors\fi\fi\fi}}\def\@evenhead{\@oddhead}
\def\@oddfoot{\small\lfoot\ifnum\count0=\startpage\copyright\ \gtp\hfill\else
\agt, Volume \thevolumenumber\ (\thevolumeyear)\hfill\fi}
\def\@evenfoot{\@oddfoot}
\makeatother
\fi
\let\maketitlepage\makeagttitle
\let\makeshorttitle\maketitlepage
\let\maketitle\maketitlepage

\newwrite\gtoutfile
\long\gdef\makeheadfile{  
{\def\\{, }\def\s{ }
\immediate\openout\gtoutfile head.xxx
\immediate\write\gtoutfile{Proxy-for: \ifx\theasciiauthors\relax
\theauthors\else\theasciiauthors\fi\s<\ifx\theasciiemail\relax\theemail\else\theasciiemail\fi>}
\immediate\write\gtoutfile{\noexpand\\}
\immediate\write\gtoutfile{Authors: \ifx\theasciiauthors\relax
\theauthors\else\theasciiauthors\fi}
{\def\\{ }\immediate\write\gtoutfile{Title: \ifx\theasciititle\relax
\thetitle\else\theasciititle\fi}}
\immediate\write\gtoutfile{Subj-class: GT or SG, GR etc}
\immediate\write\gtoutfile{MSC-class: \theprimaryclass\ifx\thesecondaryclass\relax\else, \thesecondaryclass\fi}
\immediate\write\gtoutfile{Journal-ref: Algebr. Geom. Topol. \thevolumenumber\s
(\thevolumeyear) \startpage-\finishpage}
\immediate\write\gtoutfile{Comments: Published by Algebraic and
Geometric Topology at}
\immediate\write\gtoutfile{\s\s\s  http://www.maths.warwick.ac.uk/agt/AGTVol\thevolumenumber/agt-\thevolumenumber-\thepapernumber.abs.html}
\immediate\write\gtoutfile{\noexpand\\}
\immediate\write\gtoutfile{}
\ifx\theasciiabstract\relax
\immediate\write\gtoutfile{\theabstract}\else
\immediate\write\gtoutfile{\theasciiabstract}\fi
\immediate\write\gtoutfile{}
\immediate\write\gtoutfile{\noexpand\\}
\immediate\write\gtoutfile{}
\immediate\closeout\gtoutfile}}  

\def\maketitlepage{\makeagttitle\makeheadfile}
\let\makeshorttitle\maketitlepage
\let\maketitle\maketitlepage

\lognumber{3}
\volumenumber{4}
\volumeyear{2004}
\papernumber{3}
\published{24 January 2004}
\pagenumbers{31}{47}
\received{30 May 2003}
\revised{21 August 2003}
\accepted{29 August 2003}

\usepackage{amsfonts, amsmath, graphicx}
\nocolon
\let\Bbb\mathbb

\newtheorem{pro}{Proposition}[section]
\newtheorem{thm}[pro]{Theorem}
\newtheorem{lem}[pro]{Lemma}

\newtheorem{cor}[pro]{Corollary}

\newtheorem{hyp}[pro]{Hypothesis}

\theoremstyle{definition}

\theoremstyle{remark}

\newtheorem*{qu}{Question}

\def\demo{\proof}

\def\text{\mbox}

\begin{document}

\title{Large embedded balls and Heegaard genus\\in negative curvature}
\shorttitle{Large embedded balls and Heegaard genus}
\author{David Bachman\\Daryl Cooper\\Matthew E. White}
\shortauthors{Bachman, Cooper and White}

\address{DB and MEW:\quad Mathematics Department, Cal Poly State University\\San Luis Obispo, CA 93407\\DC:\quad Mathematics Department, University of 
California\\Santa Barbara, CA 93106, USA}
\asciiaddress{DB and MEW: Mathematics Department, Cal Poly State University\\San Luis Obispo, CA 93407\\DC: Mathematics Department, University of 
California\\Santa Barbara, CA 93106, USA}

\email{dbachman@calpoly.edu, cooper@math.ucsb.edu, mewhite@calpoly.edu}

\begin{abstract} We show if $M$ is a closed, connected, orientable, 
hyperbolic 3-manifold with
Heegaard genus $g$ then $g \ge \frac{1}{2}\cosh (r)$ where $r$ denotes the
radius of any isometrically embedded ball in $M$. Assuming an 
unpublished result of Pitts and
Rubinstein improves this to $g \ge \frac{1}{2}\cosh (r) + 
\frac{1}{2}.$ We also give an upper bound on the
volume in terms of the {\em flip distance} of a Heegaard splitting, 
and describe isoperimetric surfaces in
hyperbolic balls.
\end{abstract}
\asciiabstract{%
We show if M is a closed, connected, orientable, hyperbolic 3-manifold
with Heegaard genus g then g >= 1/2 cosh(r) where r denotes the radius
of any isometrically embedded ball in M. Assuming an unpublished
result of Pitts and Rubinstein improves this to g >= 1/2 cosh(r) +
1/2.  We also give an upper bound on the volume in terms of the flip
distance of a Heegaard splitting, and describe isoperimetric surfaces
in hyperbolic balls.}

\keywords{Heegaard splitting, injectivity radius}
\primaryclass{57M50}
\secondaryclass{57M27, 57N16}
\maketitle
\let\\\par

\section{Introduction}

Let $M$ denote a closed, connected, orientable 3-manifold which 
admits a hyperbolic structure.
By Mostow rigidity the hyperbolic metric on $M$ is unique up to 
isometry. Thus geometric
invariants of the hyperbolic metric are actually topological 
invariants of $M.$ One may thus
attempt to relate geometric invariants to other topological 
invariants. For example in
\cite{cooper} it is shown that the volume of $M$ is at most $\pi$ 
times the length, $L,$ of
any presentation of the fundamental group of $M.$ A result of 
Lackenby, \cite{lack1} shows
that for alternating links the volume is bounded above and below by 
explicit affine functions of a
certain combinatorial invariant: the {\em twist number}.

The {\it injectivity radius} of $M$ is defined to be the radius of 
the smallest self-tangent
isometrically embedded ball in $M$. In \cite{white:02} White showed 
that that the injectivity
radius of $M$ is bounded above by a function of the {\it rank } of 
its fundamental group (the
minimum number of generators required to generate the group), see 
also \cite{kw} for an
extension to word hyperbolic groups. Furthermore White showed in 
\cite{white:01} that the
diameter of $M$ is bounded above, and hence the injectivity radius is 
bounded below, in terms
of $L.$

Since the Heegaard genus is always at least as large as the rank, 
White's result gives a corresponding
upper bound on injectivity radius in terms of Heegaard genus.  A 
preprint by Rubinstein
provides an improved estimate for the upper bound given by Heegaard genus
\cite{rubinstein:03}.  His techniques employ minimal surfaces in an 
intriguing way. These
ideas have inspired the result presented here.  In this note, we show 
that the radius of {\it
any} isometrically embedded ball in $M$ provides a lower bound for 
its Heegaard genus. The
precise form of our result depends on whether or not we assume an unpublished 
result (which we will refer
to as PRH) of Pitts and Rubinstein. We give two results, the sharper 
one of which
assumes PRH. In the course of the proof we construct certain 
sweepouts of manifolds of negative
curvature using simplicial surfaces whose area is bounded in terms of 
the Heegaard genus.

\begin{thm}
\label{t:main} Suppose that $M$ is a closed, orientable, connected 
Riemannian 3-manifold with
all sectional curvatures less than or equal to $-1$ and with Heegaard 
genus $g.$ Then \[g \ge
\frac{\cosh (r)}{2}.\] Assuming PRH then
\[g \ge \frac{\cosh (r) + 1}{2}.\]  Here $r$ denotes the radius of 
any isometrically embedded
ball in $M$.
\end{thm}

The proof also gives a new upper bound on the volume of a closed 
hyperbolic $3$-manifold
(\ref{volumebound}) in terms of the flip distance (see section 3) of 
a Heegaard splitting.  We have also
found it necessary to provide a statement (\ref{hypcase}) and proof 
of the fact that an equatorial disc in
a hyperbolic ball has least area among surfaces which separate the 
ball into two sets of equal volume. This
fact is certainly known to experts but does not seem to even be 
stated in the literature.

An interesting open question about 3-manifolds is the relationship 
between the Heegaard genus
and the rank of the fundamental group. Clearly
$$rank\ \le\ genus.$$ Boileau and Zieschang \cite{BZ} give an example 
of a Seifert fiber space
where the rank is $2$ but the genus is $3.$ However there is no known 
hyperbolic manifold for
which these numbers are known to be different. Furthermore, it is 
difficult to compute these
quantities. Thus the current result combined with \cite{white:02} 
provides a potential method
for constructing examples and suggests the following question:
\begin{qu} Given an integer $n$, is there a constant $R(n) > 0$ such 
that for every
closed hyperbolic 3-manifold with fundamental group generated by $n$
elements, there does not exist an isometrically embedded hyperbolic ball of
radius greater than $R(n)$?
\end{qu}

  The authors would like to thank Antonio Ros for helpful comments.  A
proper subset of the authors wish to subtitle this paper ``Big balls
imply big genus.''

The second author is partially supported by the NSF.

\section{Sweepouts}

\subsection{Heegaard Splittings} A {\it handlebody} is a 3-manifold 
which is homeomorphic to
the closure of a regular neighborhood of a graph in ${\mathbb R}^3$. 
The image of the graph
under such a homeomorphism is then a {\it spine} of the handlebody. A 
closed surface $S$ in
$M$ is a {\it Heegaard splitting} of M if $S$ separates $M$ into two 
handlebodies. We say that
the {\it genus} of a 3-manifold $M$ is the minimum among the genera 
of all Heegaard splittings
of $M$.

\subsection{Sweepouts} Let $S$ be a connected, orientable, closed 
surface and $M$ a closed, orientable
$3$-manifold. Let
$\Phi:S
\times I \rightarrow M$ be a continuous map. For each $t
\in I$ let $S_t=\Phi(S,t)$. Then $\Phi$ is a {\it sweepout} if it 
satisfies the following
properties:
\begin{enumerate}
     \item $S_0$ and $S_1$ are graphs.
     \item $\Phi_*:H_3(S\times I,\partial(S\times I))\rightarrow 
H_3(M,S_0\cup S_1)$ is an isomorphism.
\end{enumerate} We refer to condition (2) by saying that $\Phi$ is 
{\it degree one.} The
sweepout is {\em smooth} if $\Phi$ is $C^{\infty}.$  We will refer to 
the following statement
as PRH. It has been announced by Rubinstein and Pitts \cite{pr:87}, 
and a key step in the proof
has been published by Colding and De Lellis in
\cite{cold1} where they also announce that a forthcoming paper will 
complete the proof.
\begin{hyp}[PRH] Suppose that $M$ is a closed Riemannian 3-manifold 
of Heegaard genus $g.$ Then
for all $\epsilon>0$ there is a smooth sweepout $\Phi:S\times 
I\rightarrow M$ by surfaces with
$genus(S)=g$ and a minimal surface $V$ in $M$ of genus at most $g$ 
such so that $area(S_t) \le area
(V) +
\epsilon$ for all
$t\in I.$ We call such a sweepout an {\em almost minimax} sweepout.
\end{hyp}

\begin{cor}\label{PRHcor} Assume PRH. Suppose that $M$ is a closed 
orientable Riemannian
$3$-manifold of Heegaard genus $g$ such that all sectional curvatures 
are at most $-1.$ Then for every
$\epsilon>0$ there is a smooth sweepout of $M$ with the property that 
for every $t$ we have
$area(S_t)
\le 4\pi(g-1) + \epsilon.$
\end{cor}
{\bf Proof}\qua  Suppose that $V$ is a minimal surface of genus at most 
$g$ in $M.$
Since  $M$ has all sectional curvatures at most $-1$ it follows that 
$V$ has intrinsic
curvature at most $-1$.  Therefore, by the Gauss-Bonnet theorem, we have that
\[area(V) \le 2\pi(2g-2).\] Using almost minimax sweepouts provided 
by PRH gives the
result.
\qed

\medskip
We will now construct certain sweepouts which we may use instead of 
assuming PRH. The construction is
closely related to the simplicial sweepouts used by Canary and Minsky 
\cite{Can},\cite{CanMin}. In what
follows, unless otherwise stated, $M$ will always denote a closed 
Riemannian $3$-manifold with sectional
curvature at most $-1.$

We will use the term {\em triangulation} in this paper in a 
generalized sense by not requiring the
boundary of a simplex to be embedded. Formally by a {\em 
triangulation} of a manifold $M$ we mean a
CW structure on $M$ such that the corresponding
CW structure on the universal cover, $\tilde{M},$ of
$M$ is a simplicial complex and that the stabilizer (under the group 
of covering transformations) of
every simplex is trivial. By the term {\em simplex} we mean the 
closure of a cell.

We will use the terms {\em vertex, edge, triangle} and {\em 
tetrahedron} for simplices of
dimension $0,1,2$ and $3$ respectively.  In particular a {\em 
one-vertex triangulation} of a closed
orientable surface of genus $g$ consists of one
$0$-cell called the {\em vertex}, $4g-1$ one-cells whose closures are 
called {\em edges}, and whose
interiors are disjoint and both endpoints are at the vertex, and $2g$ 
two-cells whose closures are called
{\em triangles} and whose interiors are disjoint and so that the 
frontier of every triangle consists of $3$
edges.

In the constant curvature case the sweepout is by surfaces which
are a union of finitely many geodesic triangles. However we allow 
some, or all, of the triangles in a
surface to degenerate to an interval or a point. We also allow the 
triangles (or intervals) to be immersed
instead of embedded. In a general manifold of negative sectional 
curvature there are no totally geodesic
surfaces.  Instead we use triangles that are obtained by coning a 
geodesic segment to a point.

A {\em coned $n$-simplex} is a continuous map 
$\sigma:\Delta\rightarrow M$ where
$\Delta=(v_0,v_1,\cdots,v_n)$ is an $n$-simplex and $M$ is a smooth 
Riemannian manifold such that the
following holds. If $n=0$ then every such $\sigma$ is a coned 
$0$-simplex. If $n=1$ we require that
$\sigma:(v_0,v_1)\rightarrow M$ is a constant speed geodesic. In 
particular we allow the speed to be
$0$ in which case the image of $\sigma$ is a single point.  If
$n>1$ then we require:\\ (1)\qua Let $\Delta'=(v_1,\cdots,v_n)$ be the 
face of $\Delta$ which omits $v_0,$
then
$\sigma|\Delta'$ is a coned simplex.\\
(2)\qua For every $x\in\Delta'$ if we denote by $(v_0,x)$ the straight 
line in $\Delta$
regarded as a  $1$-simplex we require that $\sigma|(v_0,x)$ is a 
coned $1$-simplex.\\

We say that $\sigma$ is {\em coned from $v_0$}. Since distinct 
vertices of $\Delta$ may have the same
image under $\sigma$ we should talk about a simplex being coned from 
a {\it corner} of its image in $M.$
However, it usually causes no confusion to identify a vertex, $v_0,$ 
in $\Delta$ with its image
$\sigma(v_0)$ in $M$ and refer to the simplex as {\it coned from} 
$\sigma(v_0).$

The map $\sigma$ is smooth. We say that the simplex is
{\em degenerate} if
the image of $\sigma$ has dimension less than $n.$
Since $M$ has negative sectional curvature, every arc is homotopic keeping
  the ends fixed to a unique geodesic. Therefore, given any continuous
  $\phi:\Delta \rightarrow M$ one may homotop it, keeping the vertices fixed,
  to a unique  coned simplex. Note, however, that this map depends on the
  ordering of the vertices of $\Delta$.

In particular a coned triangle is uniquely determined by its 
restriction to the boundary together
with the choice of the vertex, $v_0,$ from which the map is coned. A 
coned triangle is a ruled
surface. Notice that we include degenerate cases: the image of a 
coned triangle might be a point or a
geodesic segment. Furthermore even if the image has dimension two, it 
may be an immersed, not embedded,
triangle. In the following discussion we will initially assume that 
all triangles are non-degenerate and
then indicate the generalization to the degenerate case.

A {\em coned simplicial surface} is continuous map $\phi:S\rightarrow 
M$ of a closed
surface, $S,$ into $M$ together with a triangulation of $S$ such that for
each triangle $T$ contained in $S$ there is
$\psi_T:(v_0,v_1,v_2)\rightarrow T$ and a coned $2$-simplex
$\sigma_T:(v_0,v_1,v_2)\rightarrow M$ such that
$\sigma_T=\phi\circ\psi_T.$ If $T$ is an embedded $2$-simplex in $S$ this
is equivalent to saying that
$\phi|T$ is a coned $2$-simplex.

The extrinsic curvature of a ruled surface is non-positive. Since all 
sectional curvatures of
$M$ are at most $-1$ it follows that the curvature of the induced 
metric on a ruled surface in $M$ also has
curvature at most $-1.$ The Gauss-Bonnet theorem and the fact that 
the sides of a
coned triangle are geodesics implies that the area of a coned 
triangle in $M$ is at most $\pi.$
If we endow $S$ with the induced (pull-back) metric then $area(S) \le 
n\pi$ where
$n$ is the number of triangles in the triangulation of $S.$

Given a vertex $v$ of the triangulation of $S$ the {\em angle sum} 
$\theta(v)$ at this vertex is the sum
of the angles of the corners of the triangles incident to $v.$ The 
Gauss-Bonnet theorem implies
$$\int_{S}\ K\ dA\ +\ \Sigma_v\ (2\pi-\theta(v))\ =\ 2\pi\chi(S).$$ 
The integral is over all of $S$ except
the edges and vertices. On this subset the Gauss curvature of the 
surface $S$ is $K\le-1.$ The sum is
over all vertices of this triangulation.
It follows that if a triangulation has at most one vertex for which the
angle sum is less than $2\pi$ then $area(S)\le \pi(4g-2).$

We can extend these ideas to the case that some, or all, of the 
triangles are degenerate by defining the
area of a degenerate triangle to be zero.

A {\em coned simplicial sweepout} is a sweepout $\Phi:S\times 
I\rightarrow M$ such that for every $t\in I$
the map $\Phi|S\times t$ is a coned simplicial surface.  The 
following result is perhaps of independent
interest.

\begin{thm}
\label{sweepout}  Suppose that $M$ is a closed orientable Riemannian
$3$-manifold of Heegaard genus
$g$ such that all sectional curvatures are at most
$-1.$ Then there is a coned simplicial sweepout $\Phi:S\times 
I\rightarrow M$ such that for every $t\in I$ the
coned simplicial surface $\Phi|(S\times t)$  consists of at most
$4g$ triangles and has at most one vertex at which the angle sum is 
less than $2\pi.$ Thus the area of
every surface in this sweepout is  at most
$\pi(4g-2).$
\end{thm}
{\bf Proof}\qua Let $S$ be a Heegaard surface of genus $g$ and $H_0,H_1$ 
the two handlebodies of the Heegaard
splitting. For $i=0,1$ choose a spine, $S_i,$ which is a wedge of $g$ 
circles in the interior of $H_i.$ We
obtain a sweepout
$\Psi:S\times I\rightarrow M$ in which $S_0$ and $S_1$ are these 
spines and for all $t \in (0,1)$ the
surfaces $S_t$ are isotopic to $S$.  Choose basepoints $p\in S$ and 
$q \in M.$ We will homotop this
sweepout to give the required coned simplicial sweepout.

If $\Psi$ is homotoped to a map $\Psi'$ so that, at each stage of the 
homotopy the image of
$S\times\partial I$ is a graph, then the degree does not change, so $\Psi'$ has
degree one and is thus a sweepout. Unless otherwise stated all the 
homotopies are of this form.

{\bf Claim 1}\qua We can homotop $\Psi$ so that $\Psi(p\times I)=q$ and 
the following holds.  For $i=0,1,$
there is a one-vertex triangulation, ${\mathcal T}_i,$ of $S\times i$ 
with vertex $p\times i$ such that
$\Psi|S\times i$ is a coned simplicial surface and each triangle of
${\mathcal T}_i$ maps to a (possibly immersed, possibly degenerate) 
geodesic arc in $M$ which starts and
ends at $q.$

\begin{figure}[ht!]
\cl{\includegraphics[width=4in]{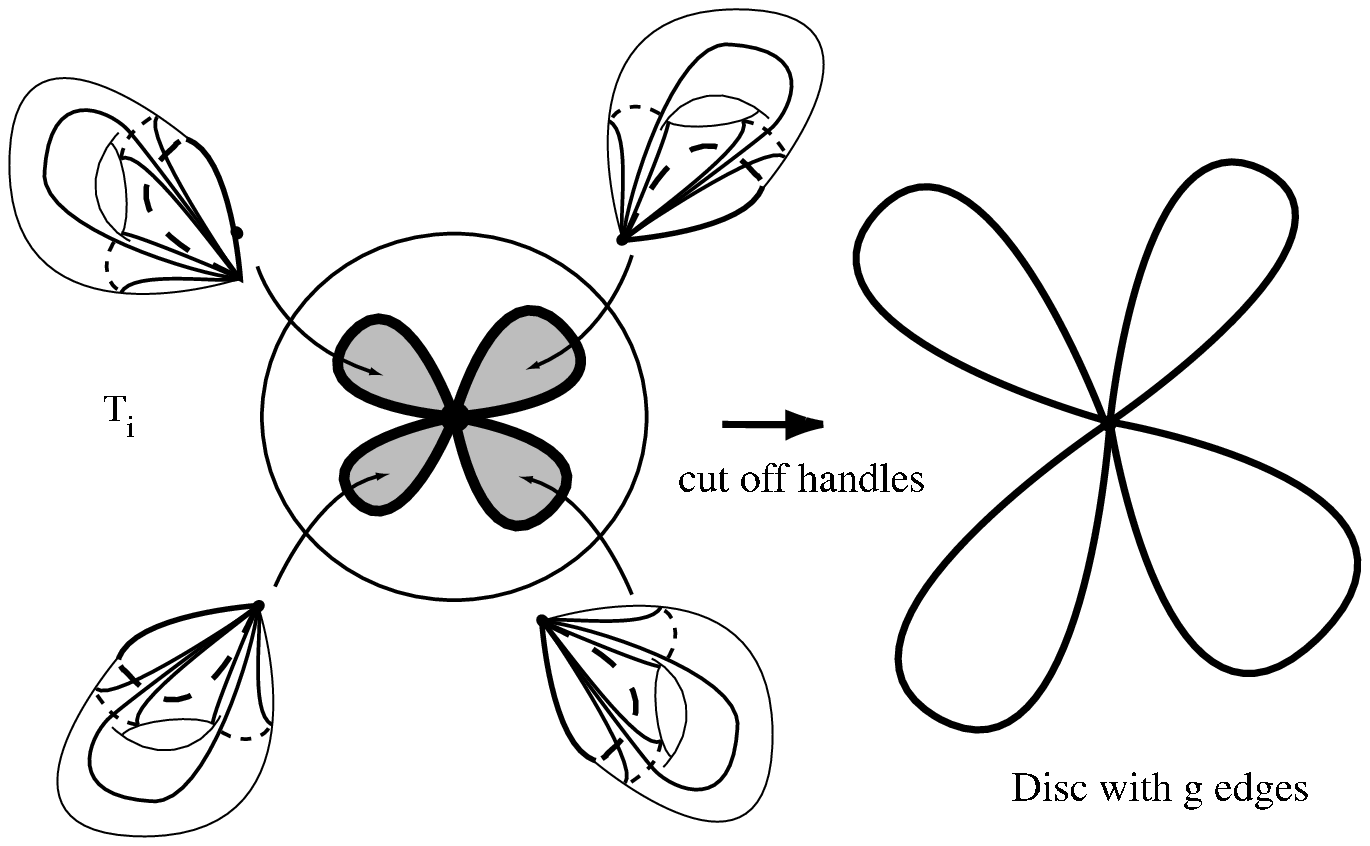}}
\caption{}
\end{figure}

We may regard ${\mathcal T}_i$ as a triangulation of $\partial H_i.$ 
The first step is to choose ${\mathcal
T}_i$ so that there is a retraction $H_i\rightarrow S_i$ so that each 
triangle of ${\mathcal T}_i$ is
mapped to either the vertex of $S_i$ or to an edge of $S_i.$  The 
proof that this can be done follows from
considering the case of genus $4$ shown in Figure 1. Each handle of 
$H_i$ contributes $4$ edges. The result
of cutting off the handles is a disc with $g$ edges. Adding an 
additional $g-3$ edges decomposes this disc
into $g-2$ triangles. [For genus $2$ a slight variation of this is used.]

First we homotop $\Psi$ so that for $i=0,1$ the restriction 
$\Psi|S\times i$ is this retraction.
Next we homotop $\Psi$ so that $p\times I$ is mapped to the basepoint 
$q$ and so that $S_i$
  is homotoped to a graph $E_i$ in $M$ which is a wedge of geodesic 
arcs based at $q.$ We
allow the case that some of these arcs are degenerate with length 
zero, and that the arcs are not
imbedded, only immersed, in $M.$ Now we may homotop $\Psi$ as a map 
of pairs $(S\times I,S\times\partial
I)\rightarrow (M,E_1\cup E_2)$ so that $\Psi|S\times i$ is a coned 
simplicial map with respect to the
triangulation
${\mathcal T}_i$ for $i=0,1.$ Each triangle of ${\mathcal T}_i$ is 
mapped to either the basepoint, $q,$ or
to an edge of $E_i.$ We denote this new sweepout by
$\Psi.$  This proves claim 1.

\medskip
{\bf Claim 2}\qua There is a map $\Psi':S\times I\rightarrow M$ which 
equals $\Psi$ on the subspace
$(S\times\partial I) \cup (p\times I)$ and such that $\Psi'$ has the 
property that for every $t\in I$ that
$\Psi'|S\times t$ is a coned simplicial surface which consists of at most
$4g$ triangles and has at most one vertex at which the angle sum is 
less than $2\pi.$

\medskip
Assuming this it follows from (\ref{homotopylemma}) (and the fact 
that $M$ is aspherical by the
Cartan-Hadamard theorem) that
$\Psi'$ is homotopic to
$\Psi$ by a homotopy fixed on $S\times\partial I.$ Thus $\Psi'$ also 
has degree one and is therefore the
required sweepout. We now prove claim 2.

\medskip
It is well known, \cite{hatcher}, that any two triangulations of a 
closed surface with one vertex
are equivalent under isotopy and a finite sequence of {\em flips.} A 
{\em flip} is
the following operation. Suppose $e$ is an edge of a triangulation 
of a surface. Then $e$
is contained in two triangles and forms the diagonal of the square 
formed by these triangles.
Remove this diagonal and replace it with the other diagonal. This 
makes sense even when some boundary
edges of the square are  identified.

\begin{figure}[ht!]
\cl{\includegraphics[width=4in]{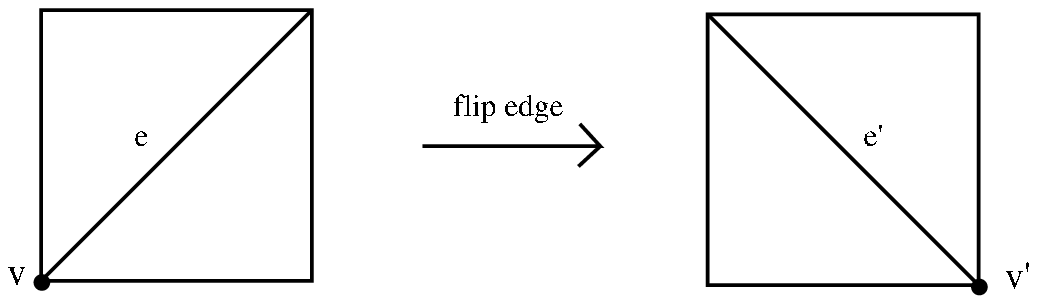}}
\caption{}
\end{figure}

Let $\{{\mathcal T}^i\}_{i=0}^n$ denote a sequence of one-vertex 
triangulations of
$S$ with the vertex at $p$ such that ${\mathcal T}^0={\mathcal T}_0$, 
${\mathcal T}^n={\mathcal T}_1$, and
for each
$i$ the triangulations
${\mathcal T}^i$ and ${\mathcal T}^{i+1}$ differ by a flip.
For each $0\le i\le n$ homotop $\Psi|S\times 0$ keeping $p$ fixed to 
obtain a coned simplicial
surface $\phi_i:S\rightarrow M$ with respect to the triangulation 
${\mathcal T}^i$ such that
$\phi_i(p)=q.$ We will choose $\phi_0=\Psi|S\times 0$ and 
$\phi_n=\Psi|S\times 1.$ Observe that
  $\phi_i$ is homotopic to $\phi_{i+1}$ keeping $p$ fixed. It is now clear that
claim 2 follows from combining the homotopies given by:

\medskip
\noindent{\bf Claim 3}\qua Suppose
  ${\mathcal T,T'}$ are one-vertex triangulations of $S$ each with 
vertex at $p$ and that ${\mathcal T'}$ is
obtained from ${\mathcal T}$ by flipping the edge $e$ of ${\mathcal 
T}$ to the edge $e'$ of ${\mathcal
T'}.$ Suppose that
$\phi,\phi':S\rightarrow M$ are coned simplicial surfaces with 
respect to these triangulations and that
they are homotopic keeping
$p$ fixed, and $\phi(p)=\phi'(p)=q.$  Then there is a homotopy
$H:S\times I\rightarrow M$ such that
$H|S\times 0=\phi$ and $H|S\times 1=\phi'$ and $H(p\times I)=q.$ 
Furthermore $H$ has the property
that for every
$t\in I$ that
$H|S\times t$ is a coned simplicial surface which consists of at most
$4g$ triangles and has at most one vertex at which the angle sum is 
less than $2\pi.$

\medskip
Let $K$ denote the union of all the edges except $e$ of ${\mathcal T},$
then $\phi|K=\phi'|K.$
Let $T$ be a triangle in ${\mathcal T}$ that does not contain $e.$ 
Then $T$ is also a triangle
in ${\mathcal T}'.$  If
$\phi$ and $\phi'$ both cone this triangle from the same vertex then
$\phi|T=\phi'|T.$ Otherwise there is a 1-parameter
family of coned simplicial surfaces connecting $\phi|T$ and 
$\phi'|T,$ see Figure 3. The maps in
this family introduce one extra vertex, $v_t,$ in one edge, $e,$ of 
$T$ and two extra triangles (each of
the triangles in
${\mathcal T}$ which contain
$v_t$ in their boundary are subdivided into two triangles).

Since $M$ has non-positive curvature, given an arc 
$\gamma:I\rightarrow M$ there is a
unique constant speed geodesic homotopic to this arc keeping 
endpoints fixed. Furthermore this
map changes continuously as the endpoints move. One then slides the 
point $v_t$ along the edge $e=[a,c]$
starting at $a$ and ending at $c$ coning from it during this process.
  The angle sum at
$v_t$ is at least
$2\pi$ because there are two edges incident to $v_t$ from opposite 
directions (their union is geodesic)
cf. the proof of lemma 4.2 of
\cite{Can}.

\begin{figure}[ht!]
\cl{\includegraphics[width=4in]{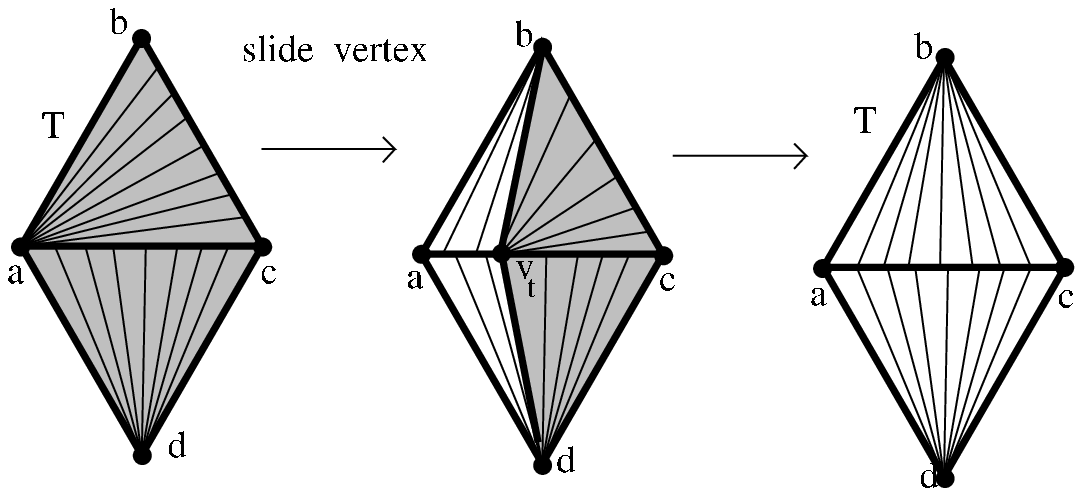}}
\caption{}
\end{figure}

There are several possible variations of the picture depending on 
which vertex each of the two triangles is
coned from. However one may view this sliding process as creating two 
new triangles (those to the left of
$v_t$ shown unshaded) and enlarging them until they entirely replace 
the two original triangles (those to
the right of
$v_t$ shown shaded.) The vertices that the new triangles are coned 
from may be chosen arbitrarily. Thus one
can use this process to change the vertices $T$ and $T'$ are coned 
from to be any chosen vertices.

Thus by doing this process at most once for each triangle we may 
arrange that $\phi$ and $\phi'$ are
equal except on the interior of the quadrilateral, $Q,$ which 
supports the flip. Let $T$ be the triangle
in ${\mathcal T}$ adjacent to $Q$ along the edge $[v,v'].$ Then $T$ 
is also a triangle in ${\mathcal T}'.$
We may slide vertices as above to arrange the coning in $Q\cup T$ is 
as shown at the left of
Figure 4: thus
$\phi|Q$ is coned from a vertex $v$ at one end of $e.$

\begin{figure}[ht!]
\cl{\includegraphics[width=3.5in]{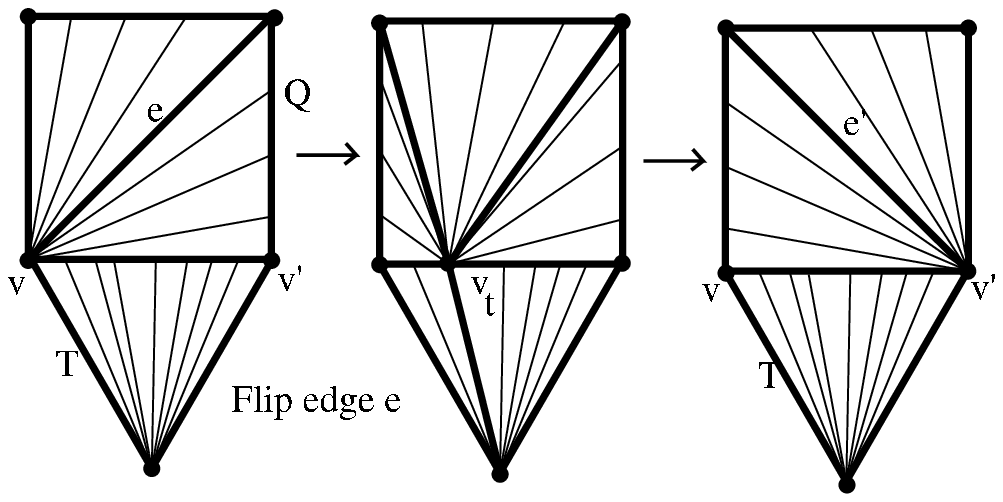}}
\caption{}
\end{figure}

  Now we can interpolate with a 1-parameter
family of coned simplicial maps which introduce one extra vertex, 
$v_t,$ and two extra triangles. As before
$\theta(v_t)\ge2\pi.$ Observe that the maps in this family do not 
change on $\partial Q.$

Finally we can use vertex sliding to change the way the triangles in 
$Q$ and $T$ are coned so that they
agree with the way they are coned by $\phi'.$ This has produced a 
sequence of coned simplicial sweepouts
interpolating between $\phi$ and $\phi'.$\qed

\begin{lem}[Homotopy lemma]\label{homotopylemma} Suppose that $X$ is 
an aspherical space and $S$
is a CW complex with a single cell of dimension $0$ called $p.$ 
Suppose that $J=[a,b]$ is an interval and
$f_0,f_1:S\times J\rightarrow X$ are continuous and that $f_0|C=f_1|C$ where
$C = (S\times \partial J)\cup (p\times J).$  Then there is a homotopy
$H:(S\times J)\times I\rightarrow X$ where
$I=[0,1]$ so that the following holds. For each $t\in I$ let 
$H_t:S\times J\rightarrow X$ be the map given
by
$H_t(s,j)=H(s,j,t)$ for $s\in S$ and $j\in J$ then $H_0=f_0$ and 
$H_1=f_1$ and $H_t|C =f_0|C$ for all
$t\in I.$
\end{lem}
\demo The hypotheses uniquely define a continuous map $H|D$ on $$D = (S\times
J\times\partial I)\cup (C\times I) = [ S\times\partial(J\times I) ] 
\cup (p\times J\times I).$$ It suffices
to show this map has a continuous extension to  $S\times J\times I.$ 
A CW structure on two spaces
determines such a structure on their product. Start with the standard 
structure on the interval with two
$0$-cells and one
$1$-cell. Use this structure on $J$ and $I$ then the product $J\times 
I$ has a structure with four
$0$-cells, four $1$-cells and one $2$-cell. We are given a cell 
structure on $S$ thus we get a product cell
structure on
$S\times J\times I$ which we now describe. For each $k$-cell $\sigma$ 
of $S$ there are:\\
(i)\qua Four $k$-cells $\sigma\times\partial J\times\partial I$ and\\
(ii)\qua Four $(k+1)$-cells $\sigma\times (\partial J\times I)$ and 
$\sigma\times( J\times \partial I)$ and\\
(iii)\qua One $(k+2)$-cell $\sigma\times J\times I.$\\

Observe that all the cells in (i) and (ii) are contained in 
$S\times\partial(J\times I)$ and give it a
CW structure. Also $p\times J\times I$ is a cell. Thus
$D$ is a subcomplex. We will extend $H$ over the cells not in $D$ in 
order of increasing dimension. Suppose
that
$\tau$ is a cell in $H$ and that $H$ has been extended over all cells 
of smaller dimension, hence $H$ is
defined on
$\partial \tau.$ Since $X$ is aspherical such an extension is always 
possible unless $dimension(\tau) =
2.$  In this case we must check that $H|\partial\tau$ is contractible 
in $X.$ However we claim that
every cell in
$S\times J\times I$ that is not contained in $D$ has dimension at 
least $3.$ It then follows that we may
extend $H$ over all cells.

Suppose that $\tau$ is a $2$-cell in $S\times J\times I$ which is not 
contained in $D.$ It follows from
the above that
$\tau =
\sigma\times J\times I$ and $\sigma$ is a $0$-cell in $S.$ This 
implies $\sigma=p,$ but $p\times J\times
I$ is contained in $D.$ This proves the claim.
\qed

\medskip
We remark here that when $M$ is hyperbolic, using lemmas 5.3 and 5.5 
of Canary \cite{Can}  one can
construct a sweepout by simplicial hyperbolic surfaces with areas at 
most $-2\pi\chi(S).$ However to extend
this construction to the case of variable negative sectional 
curvature one would have to find some
alternative to the proof of lemma 5.2 of that paper.

\medskip
\noindent{\bf Definition}\qua Suppose that $S$ is a closed, connected, 
orientable surface.
Let
$\Gamma(S)$ be the graph whose vertices are isotopy classes of 
one-vertex triangulations of $S$ and two
vertices in $\Gamma(S)$ are connected by an edge if there is an edge 
flip which converts one triangulation
to the other. Given two one-vertex triangulations of $S$ the {\em 
flip distance} between them is the least
number of edges in $\Gamma(S)$ connecting the vertices corresponding 
to these triangulations. We now use
this to define a notion of distance for a Heegaard splitting.

Suppose $M$ is a closed $3$-manifold and $S$ is a Heegaard surface 
for $M$ which separates $M$ into two
handlebodies $H_1,H_2.$ Suppose that $C_i$ is a spine for $H_i$ and 
that $C_i$ is a wedge of circles.
Suppose ${\mathcal T}_1$ and ${\mathcal T}_2$ are one-vertex 
triangulations of $S$ so that ${\mathcal T}_i$
is homotopic in $H_i$ to $C_i$ such that every edge of ${\mathcal 
T}_i$ is homotoped into either the vertex
of $C_i$ or onto one edge of $C_i.$ Then the {\em flip distance} of 
the Heegaard splitting is the
minimum over all such choices of the flip distance between ${\mathcal 
T}_1$ and ${\mathcal T}_2.$

\begin{thm}\label{volumebound} Suppose $M$ is a closed hyperbolic 
$3$-manifold and $S$ is a Heegaard
surface for $M.$  Then the hyperbolic volume of $M$ is at most
$n v_3$ where $v_3$ is the volume of a regular ideal hyperbolic 
$3$-simplex and $n$ is the flip distance of
the Heegaard splitting.
\end{thm}
{\bf Sketch Proof}\qua As in the proof of (\ref{sweepout}), given a 
sequence of edge flips one constructs a
map
$\Phi:S\times I\rightarrow M$ which equals
$\Psi$ on $S\times \partial I.$ Define a sequence of spaces $N_i$ by 
$N_0=S\times[0,1]$ and $N_{i+1}$ is
obtained from $N_i$ by glueing a tetrahedron along two of its faces 
to the quadrilateral in $\partial
N_i\setminus (S\times 0)$ which supports the $i$'th edge flip. This 
gives a cell structure for
$N_n$ relative to $S\times I$ in which the number of tetrahedra 
equals the number of edge flips.
Each tetrahedron is mapped by $\Psi$ as a coned (straight) map. Thus 
the volume of each tetrahedron is at
most
$v_3.$ The map
$\Psi$ is degree one by the homotopy lemma. This gives the volume bound.
\qed

\medskip
It is known that there are hyperbolic $3$-manifolds of Heegaard genus 
$2$ with arbitrarily large
volume; see for instance \cite{sch}.

\section{Isoperimetric surfaces}

Given a metric ball, $B,$ of radius $r$ in ${\Bbb H}^3$ an {\em 
equatorial disc} $\Delta$ is
the intersection of a hyperbolic plane with $B$ that separates $B$ 
into two halves of equal
volume.  The area of an equatorial disc is
$2\pi(\cosh(r) - 1).$

In this section we prove that any surface which divides a hyperbolic 
ball into two subsets of
equal volume has area at least as large as the area of an equatorial disc. The
literature does not seem to contain this result exactly as stated. We 
present mild variations
of proofs given by Antonio Ros \cite{ros:1},\cite{ros:2} to derive it 
from published results.
We then obtain a similar result about metrics of negative sectional 
curvature at most $-1.$

Suppose $M^n$ is a compact Riemannian manifold of dimension $n$ with 
or without boundary, and
let $V(M)$ denote its $n$-dimensional Riemannian measure which we 
will refer to as the {\em volume}. If
$\Sigma\subset M$ then the {\em area} of $\Sigma$ is the $(n-1)$-dimensional
Riemannian measure. This area will be infinite unless $\Sigma$ has 
Hausdorff dimension at most $n-1.$
The following existence and regularity results are cited as theorem 1 in
\cite{ros:2} and are the culmination of work by many authors.

\begin{thm}\label{existregular} Suppose $M^n$ is a compact Riemannian 
manifold of dimension $n$ with or
without boundary, and let $V(M)$ denote its volume. Then for any $t$ 
with $0 < t < V(M)$ there exists a
compact region $\Omega$ whose boundary $\Sigma=\partial\Omega$ 
minimizes area among regions of volume $t.$
Moreover, except for a closed singular set of Hausdorff dimension at 
most $n-8,$ the boundary, $\Sigma,$
of any minimizing region is a smooth embedded hypersurface with 
constant mean curvature, and if $\partial
M\cap\Sigma\ne\phi,$ then
$\partial M$ and $\Sigma$ meet orthogonally.
\end{thm}

  {\bf Definition}\qua Let $M$ be a Riemannian 3-dimensional manifold 
with or without boundary and
volume $V(M).$ A (possibly disconnected) surface
$\Sigma$ properly embedded in a 3-manifold
$M$ is called {\it isoperimetric} for a volume $v < V(M)$ if\\ (1)\qua 
$\Sigma$ encloses a
region of volume $v,$ and\\ (2)\qua $\Sigma$ minimizes area under the 
constraint (1).

\medskip
We shall make use of theorem 5 of
\cite{ros:1} which we quote here for the convenience of the reader:
\begin{thm}[Ros]\label{ros} {\rm(i)}\qua If the ambient space is the sphere 
$M=S^3$ with an $O(3)$-invariant
metric, then each connected component of an isoperimetric surface is 
a (topological) sphere of
revolution.\\ {\rm(ii)}\qua If
$M$ is a Euclidean ball with a radial metric, i.e. an
$O(3)$-invariant metric, then the components of any isoperimetric 
surface are either spheres
or discs of revolution.
\end{thm} A {\it surface of revolution} is a surface invariant under 
all isometric
rotations (SO(2)) about a fixed geodesic.
We apply this to the
hyperbolic metric restricted to a closed ball in hyperbolic space ${\Bbb H}^3.$
\begin{thm}\label{hypcase} Suppose that $\Omega$ is a compact subset 
of the closed ball, $B\equiv
B(r),$ of radius $r$ in ${\Bbb H}^3$ and let $S = int(B)\cap\partial 
\Omega.$ Suppose that
$Volume(\Omega) = Volume(B)/2.$  Then the area of $S$ is at least as 
large as the area of an equatorial
disc. Thus
$$area(S) \ge 2\pi(\cosh(r) - 1).$$
\end{thm} {\bf Proof}\qua
The hyperbolic metric
is homogeneous thus the existence theorem implies we may replace $S$ 
by an isoperimetric surface,
$\Sigma,$ of smaller area. Then the regularity part of the theorem 
implies that $\Sigma$ is a
smoothly embedded surface. The hyperbolic metric restricted to $B$ is 
an $O(3)$-invariant
metric on a Euclidean ball. Thus every component of
$\Sigma$ is a disc or sphere and is a surface of revolution.

We now use an argument from the proof of theorem 2 of \cite{ros:1} to 
deduce that $\Sigma$ is a
disc. First suppose that $\Sigma$ contains more than two disc 
components $\Delta$ and
$\Delta'.$ Then applying an isometry of $B$ (element of $SO(3)$) we 
may move $\Delta$
until it first becomes tangent to $\Delta'$ (or some other component 
of $\Sigma$) at a single
point. This produces a new isoperimetric surface $\Sigma'$ which is 
not regular, contradicting
the theorem. Thus
$\Sigma$ contains at most one disc component. A similar argument 
shows that if $\Sigma$
contains two components, one of which is a sphere, the sphere may be 
translated by an isometry
(of ${\Bbb H}^3$) to be tangent to some other component of $\Sigma$ 
giving the same
contradiction. Hence
$\Sigma$ is connected. It remains to show $\Sigma$ is not a sphere. 
If it were we could
translate it until it is tangent to the boundary of $B$ at a single 
point, again contradicting
regularity. Hence $\Sigma$ is a disc of revolution.

It remains to show that this disc is flat, and therefore an equatorial disc of
$B.$ Since
$\Sigma$ is a surface of revolution, $\partial\Sigma$ is a round 
circle. If this circle is the
boundary of an equatorial disc, $\Delta,$ then orthogonal projection 
onto $\Delta$ gives an
area non-increasing map $\Sigma\rightarrow\Delta,$ hence 
$area(\Sigma)\ge area(\Delta)$ which
gives the result.

Otherwise, if $\partial \Sigma$ is not the boundary of an equatorial 
disk, then since it is a round
circle, there is an equatorial disk $\Delta$ such that $\partial 
\Sigma \cap \partial \Delta
=\emptyset.$ We now use a (very slightly !) modified argument from 
the proof of theorem 5 of
\cite{ros:2}. Let $\Omega$ be the closure of the component of 
$B-\Sigma$ such that
$\partial\Omega$ is disjoint from $\partial\Delta.$ Then we may find 
a metric ball
$B'$ in ${\Bbb H}^3$ such that:\\ (1)\qquad $B'\cap\partial B = 
\Omega\cap\partial B$\\ (2)\qquad
$volume(B\cap B') = volume(\Omega) = \frac{1}{2}volume(B).$\\

 From (2) and the fact that $\Sigma$ is isoperimetric it follows that
$$area(\Sigma) \le area(B\cap\partial B').$$ Define $W=\Omega\cup 
(B'-B)\subset{\Bbb H}^3$
then $$area(\partial W) = area(\Sigma)+area(\partial B'- B) \le 
area(\partial B').$$ Now
$$volume(W) = volume(\Omega) + volume(B'-B) = volume(B').$$

\begin{figure}[ht!]
\cl{\includegraphics[width=2.5in]{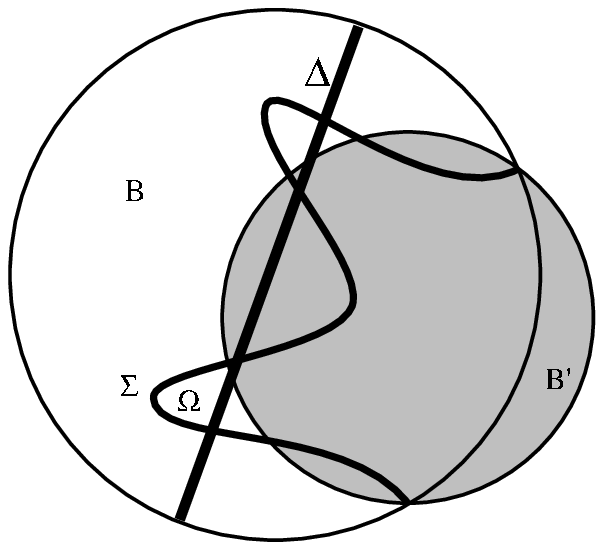}}
\caption{}
\end{figure}

It is an immediate consequence of (\ref{existregular}) and 
(\ref{ros}) that given $v>0,$ if a closed
surface in
${\Bbb H}^3$ bounds a region of volume $v$ and it has least area 
among such surfaces then it is a
hyperbolic sphere. It follows that
$W$ is a metric ball in
${\Bbb H}^3$ and thus $W=B'.$ Hence $\Sigma = B\cap \partial B'$ is a 
subset of the hyperbolic sphere
$\partial B'.$

Since $\Sigma$ is orthogonal to $\partial B$ it follows that 
$\partial B'$ is orthogonal to $\partial
B.$ But then $\partial B'$ must be disjoint from $\Delta$ and $B'\cap 
B$ is strictly contained in one
component of $B-\Delta.$ This contradicts that $B'\cap B$ has the 
same volume as that component.\qed

\begin{cor}\label{varycase}Suppose that $N$ is a complete simply 
connected Riemannian
3-manifold with all sectional curvatures everywhere less than or 
equal to $-1.$ Let $B\equiv
B_r(p;N)$ be a ball of radius $r>0$ centered at some point $p\in N.$ 
Suppose $q\in {\Bbb H}^3$ and let
$B'=B_r(q;{\Bbb H}^3)$ denote the ball of radius $r$ center $q$ in 
${\Bbb H}^3.$  Let
$\phi:T_pN\rightarrow T_q{\Bbb H}^3$ be an isometric linear map. Then 
the exponential maps
$$exp_p:T_pN\rightarrow N\qquad\qquad exp_q:T_q{\Bbb 
H}^3\rightarrow{\Bbb H}^3$$ are diffeomorphisms and
$\phi(B)=B'.$
Suppose that $\Omega\!\subset\! B$ and that
$volume(\phi\Omega)$ $= Volume(B')/2.$ Define $S=\partial\Omega\setminus\partial
B$ then  $area(S)\ge2\pi(\cosh(r) - 1).$
\end{cor}

{\bf Proof}\qua  By the Cartan-Hadamard theorem the exponential maps are 
diffeomorphisms. Define
$g=exp_q\circ\phi\circ exp_p^{-1}:N\rightarrow {\Bbb H}^3.$  Since
$\phi(0)=0$ it follows that $g(p)=q$ and since the exponential maps 
are isometries on one-dimensional
subspaces if follows that $\phi(B)=B'.$ Since
$K\le-1$ the map $g$ is distance non-increasing. Thus it is area 
non-increasing, hence
$area(S)\ge area(g(S)).$ By the theorem $area(g(S))\ge2\pi(\cosh(r) - 1).$
\qed

\medskip

This suggests the following: \begin{qu} Suppose that $B$ is a closed 
ball of radius $r$ in a
complete simply connected Riemannian 3-manifold with all sectional 
curvatures less than or
equal to $-1.$ Suppose that $\Sigma$ is an isoperimetric surface in 
$B$ which separates it into
two open sets of equal volume. Is the area of $\Sigma$ at least 
$2\pi(\cosh(r) - 1)$?
\end{qu}
\section{Proof of Theorem (\ref{t:main})}

\begin{thm}\label{areabound} Suppose that $M$ is a closed, connected, 
orientable Riemannian $3$-manifold with
all sectional curvatures at most $-1.$ Suppose that there is a 
piecewise smooth sweepout of
$M$ by surfaces all of which have area at most $A.$ Suppose that 
$p\in M$ and $r=inj(p),$ then
$$2\pi(cosh(r)-1) \le A.$$
\end{thm}
{\bf Proof}\qua We first consider the case that $M$ is hyperbolic. Let 
$B$ denote an
isometrically embedded, standard ball of radius $r$ in $M$. We claim 
that there is a (possibly
non-embedded) surface of area at most $A$ which separates $B$ into 
two open sets of equal volume. It then
follows from (\ref{hypcase}) that $2\pi(cosh(r)-1) \le A.$

Let $\Phi:S\times I\rightarrow M$ be the given sweepout.
Initially suppose that $\Phi|S\times(0,1)$ is an embedding.
Consider the function
$v:I\rightarrow {\Bbb R}$ defined by $v(t)=Volume(B\cap S[0,t]).$
Then $v$ is continuous and $v(0)=0.$ Since a sweepout has
degree $1$ it follows that
$v(1) =  volume(B).$ Hence there is some $t=t_0$ such that
$v(t_0)=\frac{1}{2}volume(B).$ Then $B\cap\Phi(S\times t_0)$ is the 
required surface.

The general hyperbolic case is handled as follows. The compact set
$S_t=\Phi(S\times t)$ has finite area and therefore volume $0.$ Let
$U_t$ be a component of
$M\setminus S_t.$ Then $U_t$ is open and its frontier is contained in 
$S_t.$ Consider
$V_t=(S\times[0,t))\cap\Phi^{-1}(U_t).$ Then $\Phi|:V_t\rightarrow 
U_t$ is proper and thus has a well
defined degree. Let $\Omega(t)$ denote the closure of the union of 
those components $U_t$ for which this
degree equals
$1.$ Then  $v(t)=Volume(B\cap \Omega(t))$ is continuous in $t.$ When 
$t=0$ clearly $\Omega(0)=\phi.$ When
$t=1,$ since a sweepout has degree one $\Omega(1)=M.$ Thus
$volume(B\cap\Omega(1))=Volume(B).$ Thus there is
$t$ with
$v(t)=\frac{1}{2}volume(B).$ The frontier of $\Omega(t)$ is contained 
in $S_t$ and is thus a piecewise smooth
surface with area at most $A.$ This surface might not be embedded in 
$M,$ however theorem (\ref{hypcase})
still applies (it
  does not assume that the boundary of the region is an embedded 
surface.) This completes the proof in the
case that
$M$ is hyperbolic.

In the variable curvature case, we first use the map $g$ from the 
proof of (\ref{varycase}) to
use the given sweepout of $M$ to construct a sweepout of a ball in 
${\Bbb H}^3$ and then use
(\ref{varycase}) to obtain the result.
\qed

\medskip
The proof of Theorem (\ref{t:main}) is now completed by noting that 
the first inequality follows directly
from Theorems (\ref{areabound}) and (\ref{sweepout}), while the 
second follows from a combination of
Theorem (\ref{areabound}) and Corollary (\ref{PRHcor}).

\Addresses\recd

\end{document}